\documentclass[a4paper,12pt,twoside,reqno]{amsart}

\usepackage[pagebackref]{hyperref}
\usepackage{oitp, mymacros}

\title[Calabi--Yau complete intersections]{Calabi--Yau complete intersections \\
in exceptional Grassmannians}

\author[A.~Ito]{Atsushi Ito}
\address{Department of Mathematics,
Faculty of Environmental, Life, Natural Science and Technology,
Okayama University,
1-1-1 Tsushima-naka, Kita-ku, Okayama-shi,
700-8530, Japan}
\email{ito-atsushi@okayama-u.ac.jp}

\author[M.~Miura]{Makoto Miura}
\address{
Department of Mathematics,
Graduate School of Science,
Osaka University,
Machikaneyama 1-1,
Toyonaka,
Osaka,
560-0043,
Japan}
\email{miurror.jp@gmail.com}

\author[S.~Okawa]{Shinnosuke Okawa}
\address{
Department of Mathematics,
Graduate School of Science,
Osaka University,
Machikaneyama 1-1,
Toyonaka,
Osaka,
560-0043,
Japan}
\email{okawa@math.sci.osaka-u.ac.jp}

\author[K.~Ueda]{Kazushi Ueda}
\address{
Graduate School of Mathematical Sciences,
The University of Tokyo,
3-8-1 Komaba,
Meguro-ku,
Tokyo,
153-8914,
Japan}
\email{kazushi@ms.u-tokyo.ac.jp}

\dedicatory{Dedicated
with admiration
to the
memory of
Bumsig Kim}

\begin{document}

\begin{abstract}
We classify
completely reducible 
equivariant vector bundles
on Grassmannians of exceptional Lie groups
which give Calabi--Yau 3-folds
as complete intersections.
We also calculate Hodge numbers for those Calabi--Yau 3-folds.
\end{abstract}

\maketitle
\section{Introduction}

A smooth projective manifold
is said to be \emph{Calabi--Yau}
if the canonical bundle 
is trivial.
Calabi--Yau manifolds
have attracted attentions
from both mathematicians and string theorists,
not only because of their importance in the classification of algebraic varieties,
but also because of their relation with string theory and mirror symmetry.
A Calabi--Yau manifold in dimensions at most two are
either an elliptic curve,
an abelian surface, or a K3 surface.
In dimensions greater than two,
it is not known whether the number of deformation equivalence classes
(or even homeomorphism types)
of Calabi--Yau manifolds is finite or not.


In this paper,
we study Calabi--Yau 3-folds
in rational homogeneous spaces
of exceptional types.
The main result is the following:

\begin{theorem} \label{th:classification}
A complete intersection Calabi--Yau 3-fold 
of a globally generated completely reducible 
equivariant vector bundle $\cE$
on an exceptional Grassmannian $G/P$,
which is not 
a complete intersection of line bundles
on a projective space,
is one of those appearing in \pref{tb:CY3}.
\end{theorem}

\begin{table}[h]
\centering
\begin{align*}
\begin{array}{ccccccc}
\toprule
\text{No.} & G/P & \cE &h^{1,1}&h^{1,2}&\deg&c_2\\ 
\midrule
1 & \dynkin[parabolic=1]G2 
& (1,1) 
&1&50&42&84\\ 
2 & \dynkin[parabolic=2]G2 
& (1,1) 
&1&50&14&56\\ 
3& \dynkin[parabolic=12356]E6 
& (1,0,0,0,0;0)\oplus(0,0,0,0,1;0)^{\oplus 4}
&1&31&192&132 \\
4 & \dynkin[parabolic=1]G2 
& (1,0)\oplus (2,0) 
&1&61&36&84\\
\bottomrule 
\end{array}
\end{align*}
\caption{Complete intersection Calabi--Yau 3-folds in exceptional Grassmannians}
 \label{tb:CY3}
\end{table}

In particular,
there is no such Calabi--Yau 3-fold
in exceptional Grassmannians
of types $E_7$, $E_8$, and $F_4$.

In \pref{tb:CY3},
we label the simple roots of $E_6$ as
\[
\dynkin[labels={\alpha_1,\alpha_6,\alpha_2,\alpha_3,\alpha_4,\alpha_5},edge length=.75cm]E6
\]
and write the coordinates
of a weight
$
\lambda = \lambda_1 \omega_1 + \cdots + \lambda_6 \omega_6
$
with respect to the corresponding fundamental weights
$\omega_1, \ldots, \omega_6$
as
$
(\lambda_1,\ldots,\lambda_5;\lambda_6).
$
Similarly,
we label the simple roots as
\[
\dynkin[labels={\alpha_1,\alpha_2},edge length=.75cm]G2
\]
and write
$
\lambda
= \lambda_1 \omega_1 + \lambda_2 \omega_2
= (\lambda_1, \lambda_2)
$
for $G_2$.
A weight $\lambda$ is identified
with the equivariant vector bundle
associated with the irreducible representation
of the Levi subgroup
with highest weight $\lambda$.
The degree and the second Chern number are
with respect to the restriction of the ample generator
of the Picard group of the ambient space.


Note that the classification itself in 
\pref{th:classification} was already obtained in \cite{benedettiphd}.
The Hodge numbers are newly calculated in this work, and are also reverified in \cite[Table 2]{2108.13314} shortly after that.

The $G_2$-Grassmannian
\dynkin[parabolic=1]G2
is the zero locus of the section
$s \in H^0(\cQ^\dual(1)) \cong \bigwedge^3 \bC^7$
corresponding to the $G_2$-invariant 3-form,
and $\cE_{(1,1)}$ is the restriction of $\cS^\dual(1)$,
where $\cS$ and $\cQ$ are the universal subbundle
and the universal quotient bundle on $\Gr(2,7)$.
Hence No.~1 in \pref{tb:CY3} is the same as No.~16
in \cite[Table 1]{MR3968921},
which is known
by \cite[Proposition 5.1]{MR3968921}
to be deformation-equivalent
to the intersection of the image of $\Gr(2,7)$
and a linear subspace of codimension 7
in $\bP \lb \bigwedge^2 \bC^7 \rb$.

The $G_2$-Grassmannian
\dynkin[parabolic=2]G2
is a smooth quadric hypersurface in $\bP^6$.
Calabi--Yau 3-folds contained in a (not necessarily smooth) quadric
5-fold are classified in \cite[Section 5]{MR3476687},
and it is shown in \cite[Theorem 7.1]{MR3476687} that
No.~2 in \pref{tb:CY3}
is deformation-equivalent
to the Pfaffian Calabi--Yau 3-fold
appearing in \cite{MR1775415}.

Although the families
No.~1 and No.~2 were known,
their relation with $G_2$-Grassmannians were new,
and led to the discoveries of an L-equivalence
\cite{MR3912058},
a derived equivalence
\cite{MR3830796},
and a 7-fold flop
\cite{MR3959275}.

To the best of our knowledge,
no known Calabi--Yau 3-fold
has the same topological invariants as 
No.~3 in \pref{tb:CY3}.
The restriction $\cO_X(1)$ of the ample generator of the ambient space is primitive
since
\begin{align}
  \chi(\cO_X(t)) = 32 t^3 + 11 t
\end{align}
by the Hirzebruch--Riemann--Roch theorem \pref{eq:HRR}.

Families of Calabi--Yau 3-folds 
described as complete intersections of line bundles 
on the $G_2$-Grassmannian
\dynkin[parabolic=2]G2
are omitted in \pref{th:classification},
since they are complete intersections of line bundles in $\bP^6$.

Calabi--Yau complete intersection 3-folds
of completely reducible
equivariant bundles
on the Cayley plane
\dynkin[parabolic=1]E6
containing $\scE_{\omega_5}$ as a direct summand
turn out to be the empty set,
since the zero locus of a general section
of $\scE_{\omega_5}$
on the Cayley plane
is the empty set because $h^0(\cO_X)=0$.

The proof of \pref{th:classification} also shows the following:

\begin{theorem} \label{th:Fano}
  There is no complete intersection Fano 4-fold 
  of a globally generated completely reducible 
  equivariant vector bundle
  on an exceptional Grassmannian,
  which is neither
  a complete intersection of line bundles
  on a projective space
  nor a hypersurface of the $G_2$-Grassmannian \dynkin[parabolic=1]G2.
\end{theorem}

We also classify a class of Calabi--Yau 3-folds
in flag varieties of exceptional types.

\begin{theorem} \label{th:classification2}
A globally generated
completely reducible
equivariant vector bundle $\cE$
on an exceptional flag variety $G/P$
of Picard number greater than one
satisfying $\rank \cE = \dim G/P -3$ and $c_1(\cE) = c_1(G/P)$
is one of those
appearing in \pref{tb:CY3flag}.
\end{theorem}

\begin{table}
  \centering
  \begin{align*}
  \begin{array}{ccccc}
  \toprule
  \text{No.} & G/P & \cE &h^{1,1}&h^{1,2}\\
  \midrule
  5 & \dynkin[mark=x]G2 
  & (1,0)^{\oplus 2}\oplus (0,2) 
  &2&38\\ 
  6 & \dynkin[mark=x]G2 
  & (1,0)\oplus (0,1)\oplus (1,1) 
  &2&48\\ 
  7 & \dynkin[mark=x]G2 
  & (2,0)\oplus (0,1)^{\oplus 2}
  &2&58\\
  \bottomrule 
  \end{array}
  \end{align*}
  \caption{Complete intersection Calabi--Yau 3-folds in exceptional flag manifolds}
  \label{tb:CY3flag}
\end{table}

The condition $c_1(\cE) = c_1(G/P)$ is sufficient
for the zero locus $X$ to be Calabi--Yau.
This condition will be necessary
if the restriction $\Pic G/P \to \Pic X$ is injective,
and it is an interesting problem to decide when this is the case.

\emph{Acknowledgements.}
We thank Daisuke Inoue for helpful discussions on Pfaffian Calabi--Yau 3-folds,
Grzegorz Kapustka and Michał Kapustka for pointing out the reference
\cite{MR3476687},
Pieter Belmans and Maxim Smirnov for correcting an error 
about representation theory of parabolic subalgebra and 
for informing us the reference \cite{MR124919},
and Kyeong-Dong Park for pointing out that $\cE_{\omega_5}$
defines the empty set in the Cayley plane.
A.~I.~was supported by Grant-in-Aid for Scientific Research
(14J01881,
17K14162).
M.~M.~was supported by
Frontiers of Mathematical Sciences and Physics at University of Tokyo,
Korea Institute for Advanced Study, and Grant-in-Aid for Scientific Research
(21K03156).
S.~O.~was partially supported by Grants-in-Aid for Scientific Research
(16H05994,
16K13746,
16H02141,
16K13743,
16K13755,
16H06337)
and the Inamori Foundation.
K.~U.~was partially supported by Grants-in-Aid for Scientific Research
(24740043,
15KT0105,
16K13743,
16H03930).

\section{Equivariant vector bundles over \texorpdfstring{$G/P$}{G/P}}
 \label{sc:G/P}
Let $\frakg$ be a complex semisimple Lie algebra of rank $r$.
The corresponding simply-connected Lie group is denoted by $G$.
Fix a Cartan subgroup $H \subset G$
with the associated Cartan subalgebra $\frakh \subset \frakg$ and set
\begin{align}
 \frakg_\alpha \coloneqq \lc v \in \frakg \relmid [h,v]=\alpha(h)v
  \text{ for any } h \in \frakh \rc
\end{align}
for
$
 \alpha \in \frakh^\dual \coloneqq \Hom_\bC(\frakh,\bC).
$
One has the root decomposition
\begin{align}
 \frakg = \frakh \oplus \bigoplus_{\alpha\in \Delta} \frakg_\alpha, 
\end{align}
where
\begin{align}
 \Delta \coloneqq
  \lc \alpha \in \frakh^\dual \relmid \frakg_\alpha \ne \lc 0 \rc, \alpha \ne 0 \rc.
\end{align}
We choose a system of simple roots
$
 \cS \coloneqq \lc \alpha_1, \dots, \alpha_r \rc \subset \Delta.
$
This choice is equivalent to the choice of the sets
$\Delta^+$ and $\Delta^-$
of positive and negative roots.

The Dynkin diagram of $\frakg$ is a graph
whose nodes correspond to
the simple roots $\alpha_i \in \cS$
and whose edges represent the Cartan integers
$\killing{\alpha_i}{\alpha_j^\dual}$,
where
$\killing{-}{-}$ is the Killing form on $\frakh^\dual$ and
$
 \alpha^\dual = 2 \alpha / \killing{\alpha}{\alpha}$.
One has
$
 \killing{\alpha_i}{\alpha_i^\dual} = 2
$
for all $\alpha_i \in \cS$, and
the correspondence between edges 
and the Cartan integers is given by
\begin{align}
&\dynkin[labels*={\alpha},label height=d,label depth=d]A1 \, 
 \dynkin[labels*={\beta},label height=d,label depth=d]A1 
\quad \Longleftrightarrow \quad 
\killing{\alpha}{\beta^\dual} =
\killing{\beta}{\alpha^\dual} = 0, \\
&\dynkin[labels*={\alpha,\beta},label height=d,label depth=d]A2
\quad \Longleftrightarrow \quad
\killing{\alpha}{\beta^\dual} =
\killing{\beta}{\alpha^\dual} = -1, \\
&\dynkin[labels*={\alpha,\beta},label height=d,label depth=d]B2
\quad \Longleftrightarrow \quad 
\killing{\alpha}{\beta^\dual} = -2, \quad
\killing{\beta}{\alpha^\dual} = -1, \\
&\dynkin[labels*={\alpha,\beta},label height=d,label depth=d]G2
\quad \Longleftrightarrow \quad 
\killing{\alpha}{\beta^\dual} = -3, \quad
\killing{\beta}{\alpha^\dual} = -1.
\end{align}

A subgroup $P$ of $G$ is said to be \emph{parabolic}
if $G/P$ is a projective variety.
Conjugacy classes of parabolic subgroups
are in one-to-one correspondence
with subsets $\cS_\frakp \subset \cS$ of the set of simple roots
in such a way that the corresponding subalgebra $\frakp$ is given by
\begin{align}
 \frakp \coloneqq \frakl \oplus \frakn,
\end{align}
where the \emph{Levi part} $\frakl$ is
\begin{align}
 \frakl \coloneqq \frakh \oplus
  \bigoplus_{\alpha \in \lb \vspan \cS_\frakp \rb \cap \Delta} \frakg_\alpha,
\end{align}
and the \emph{nilpotent part} $\frakn$ is
\begin{align}
 \frakn \coloneqq \bigoplus_{\alpha \in \Delta^+ \setminus \vspan \cS_\frakp} \frakg_\alpha.
\end{align}
Here, $\vspan \cS_\frakp \subset \frakh^\dual$ is the linear subspace
spanned by $\cS_\frakp$.
The subset $\cS_\frakp \subset \cS$ can be described
by a \emph{crossed Dynkin diagram},
where elements not in $\cS_\frakp$ are crossed out
(i.e., elements of $\cS_\frakp$ correspond to uncrossed nodes).
The inclusion relation of $\cS_\frakp$ corresponds to the inclusion relation of $P$.
For example,
the Borel subgroup is the minimal parabolic subgroup,
so that all the nodes are crossed out
in the corresponding crossed Dynkin diagram.
We write the Weyl group of $\frakl$ as $W_P$,
which is the subgroup of $W=W_G$
generated by simple reflections associated with elements of $\cS_\frakp$.
One has
\begin{align}
 \dim G/P
  = \# \left( \Delta^- \setminus \vspan \cS_\frakp\right)
  = \dim G/B - \dim G'/B', 
\end{align}
where $G'/B'$ is the full flag variety
corresponding to the full subgraph
of the Dynkin diagram of $G$
consisting of uncrossed nodes.

A finite-dimensional representation 
$V$ of the parabolic subalgebra $\frakp$
naturally carries a filtration
\begin{align}
  V = F^1V \supset \cdots \supset  F^jV \supset F^{j+1}V
  \supset \cdots \supset F^s V
  \supset F^{s+1}
  = F^{s+2} = \cdots
  = 0
\end{align}
for some $s$,
where we set
\begin{align}
  F^jV = \frakn^{j-1}\cdot V
\end{align}
for all $j$. The nilpotent part $\frakn$ acts trivially 
on each quotient
\begin{align}
  V_j \coloneqq F^jV / F^{j+1}V,
\end{align}
so that one may regard $V_j$
as a representation of the Levi subalgebra $\frakl$.
Since $\frakl$ is reductive,
any finite-dimensional representation of $\frakl$ is completely reducible,
i.e.,
the direct sum of irreducible representations.
We write
\begin{align} \label{eq:sum}
  V
  = \sum_j V_j
  = V_1 + V_2 + \cdots + V_s,
\end{align}
where
$
F^j V
= V_j + \cdots + V_s
$
is a sub-representation.
Let $W = \sum_j W_j$ be another finite-dimensional representation of $\frakp$.
It is easy to see
\begin{align}
  V \oplus W
  &= \sum_{j\ge 1} \lb V_j\oplus W_j \rb,\\
  V \otimes W
  &= \sum_{l\ge 2} \bigoplus_{j+k=l} \lb V_j\otimes W_k \rb,\\
  \label{eq:exterior}
  \bigwedge^p V
  &= \sum_{k\ge p} \bigoplus_{h(\bfp)=k} 
  \bigotimes_{j\ge 1} \bigwedge^{p_j}V_j,\\ 
  \label{eq:symmetric}
  \Sym^p V
  &= \sum_{k\ge p} \bigoplus_{h(\bfp)=k} 
  \bigotimes_{j\ge 1} \Sym^{p_j}V_j,
\end{align}
for any positive integer $p$,
where $\bfp=(p_1,\dots, p_r)$
runs over partitions of $p$ and
$h(\bfp)= \sum_j jp_j$.

Irreducible representations of the parabolic subalgebra $\frakp$
correspond bijectively to those of the Levi subalgebra $\frakl$,
which are highest weight representations.
Since $\frakg$ and $\frakl$ share the same Cartan subalgebra,
weights of $\frakl$ can be regarded naturally as weights of $\frakg$.
Since $G$ and $L$ share the same Cartan subgroup,
the notion of integrality of weights is the same for both $\frakg$ and $\frakl$.
The fundamental weight associated with the simple root $\alpha_i \in \cS$
is denoted by $\omega_i$;
\begin{align}
 \killing{\omega_i}{\alpha_j^\dual} = \delta_{ij}, \quad i, j = 1, \ldots, r.
\end{align}
A weight $\lambda = \sum_{i=1}^r \lambda_i \omega_i$ is
\begin{itemize}
 \item
\emph{integral}
if $\lambda_i \in \bZ$ for any $i = 1, \ldots, r$,
 \item
\emph{$\frakp$-dominant}
if $\lambda_i \in \bN$ for any $i$ such that $\alpha_i \in \cS_\frakp$, and
 \item
\emph{$\frakg$-dominant}
if $\lambda_i \in \bN$ for any $i = 1, \ldots, r$.
\end{itemize}
A highest weight representation of $\frakl$ integrates to a representation of $P$
if and only if the highest weight is integral and $\frakp$-dominant
(see e.g. \cite[Remark 3.1.6]{MR1038279}).
The irreducible representations of $P$ and $G$
with highest weight $\lambda$
are denoted by $\VP_\lambda$ and $\VG_\lambda$ respectively.

The category of equivariant vector bundles on $G/P$ 
is equivalent to that of representations of $P$.
For a representation $V$ of $P$,
the corresponding equivariant vector bundle on $G/P$ is denoted by
\begin{align}
 \cE_V \coloneqq G \times_P V.
\end{align}
We write the filtration of an equivariant bundle
corresponding to the filtration \pref{eq:sum}
as
\begin{align}
  \label{eq:sumbundle}
  \cE = \sum_{j\ge 1} \cE_j = \cE_1 + \cE_2 + \cdots + \cE_s.
\end{align}
For an irreducible equivariant vector bundle,
we set
\begin{align}
 \cE_\lambda
  \coloneqq \cE_{\VP_\lambda}^\dual,
\end{align}
which is globally generated
if and only if $\lambda$ is $\frakg$-dominant.
The Borel--Weil--Bott theorem gives
an isomorphism
\begin{align} \label{eq:Borel--Weil--Bott}
 H^{\ell(w)} \lb G/P, \cE_\lambda \rb
  \cong \lb \VG_{w.\lambda} \rb^\dual
\end{align}
of $G$-vector spaces,
where
$
 \rho \coloneqq 1/2 \sum_{\alpha \in \Delta^+} \alpha
$
is the Weyl vector,
\begin{align}
 w.\lambda \coloneqq w(\lambda+\rho)-\rho
\end{align}
is the affine Weyl action,
$w \in W_P$ is the unique element
such that
$w.\lambda$ is $\frakg$-dominant
(the left hand side of \pref{eq:Borel--Weil--Bott} is zero
if there is no such $w$),
and $\ell(w)$ is the length of (the minimal representative in $W$ of)
$w \in W_P$.
The filtration \pref{eq:sumbundle}
gives a spectral sequence
\begin{align}
  \label{eq:extension}
  E_1^{p,q}=H^{p+q}(G/P, \cE_p) \Rightarrow H^{p+q}(G/P,\cE),
\end{align}
which allows us to compute the cohomology of $\cE$ using
the Borel--Weil--Bott theorem.

The Picard group $\Pic G/P$ is isomorphic
to the group $\Hom(P, \bCx)$ of characters of $P$.
The set of weights of a representation $V$ of $P$
is denoted by $\Delta(V)$.
One has
\begin{align}
 \rank \cE_V = \dim V = \no{\Delta(V)}, \ 
 \det \cE_V \cong \cE_{\det V}, \text{ and }
 \Delta(\det V) = \lc \sum_{\lambda \in \Delta(V)} \lambda \rc.
\end{align}
The tangent bundle $T_{G/P}$ corresponds
to the $P$-vector space
$\frakg/\frakp$
with respect to the adjoint action;
\begin{align}
 T_{G/P} \cong \cE_{\frakg/\frakp}.
\end{align}
Since the weights of the adjoint action are roots,
one has
\begin{align}
 \Delta(\frakg/\frakp)
  = \Delta(\frakg) \setminus \Delta(\frakp)
  = \Delta^- \setminus \vspan \cS_\frakp.
\end{align}
Although
the tangent bundle $T_{G/P}$
and hence the cotangent bundle $\Omega^1_{G/P}$
are indecomposable if $G$ is simple,
they are not irreducible unless
$G/P$ is a Hermitian symmetric space \cite[Theorem 6]{MR124919}.
Instead,
it carries the following filtration.
Recall that the \emph{height} of a positive root
is the sum of the coefficients of the simple roots. 

\begin{lemma}[{cf.~e.g.~\cite[Section 9.9]{MR1038279}}] \label{lem:cotangent}
  For a complex semisimple simply-connected Lie group $G$
  and
  a parabolic subgroup
  $P\subset G$,
  the cotangent bundle has the filtration
  \begin{align}
    \label{eq:decomposition}
    \Omega^1_{G/P} = \sum_{j\ge 0} \bigoplus_{\abs{\bfn}=j} \cE_{-\alpha(\bfn)},
  \end{align}
  where $\bfn=(n_\alpha)_\alpha$ runs over the image of the map
  \begin{align}
    \pi \colon \Delta \lb \lb \frakg/\frakp \rb^\vee \rb
    = \Delta^+\setminus \vspan \cS_\frakp \rightarrow 
    \bZ_{\ge 0}^{\cS\setminus \cS_\frakp}
  \end{align}
  taking the coefficients of the roots 
  in $\cS\setminus \cS_\frakp$,
  $\abs{\bfn} \coloneqq \sum_{\alpha \in \cS\setminus \cS_\frakp}n_\alpha$,
  and $\alpha(\bfn)\in \Delta^+\setminus \vspan \cS_\frakp$ is 
  the unique element of minimal height such that $\pi(\alpha)=\bfn$.
\end{lemma}

\begin{proof}[Sketch of proof]
  Since the coadjoint action of $\frakn$ on
  $
  \lb \frakg/\frakp \rb^\vee
  \subset
  \frakg^\vee
  $ 
  increases $\abs{\bfn}$,
  it suffices
  to show the decomposition for each graded component.
  For each $\bfn \in \Ima \pi$, the root $\alpha(\bfn)$ is uniquely determined,
  since
  the existence of distinct positive roots 
  $\beta_1, \beta_2$
  with the same height
  contradicts the fact that
  the difference
  $\beta_1 - \beta_2$ must be a
  (positive or negative) root
  \cite[Lemma 9.4]{MR0323842}.
  For any $\beta \in \pi^{-1}(\bfn)$,
  the difference $\alpha(\bfn) - \beta$ is a root in $\vspan \cS_\frakp$,
  and the ladder operator of the corresponding $\mathfrak{sl}_2$-triple
  sends
  $
  \bC(-\alpha(\bfn))
  \subset
  \lb \frakg/\frakp \rb^\vee
  \subset
  \frakg^\vee
  $
  onto
  $\bC (-\beta)$
  by the coadjoint action \cite[Proposition 8.4]{MR0323842}.
  This implies the existence of the filtration \pref{eq:decomposition}.
\end{proof}

\section{Complete Intersections of Equivariant Vector Bundles}
Let $\cE \coloneqq \cE_V$ be the equivariant vector bundle
on $F \coloneqq G/P$
associated with a representation $V$ of $P$.
Assume that $\cE$ is globally generated.
For a general section $s$ of $\cE$,
the zero locus $X \coloneqq s^{-1}(0)$
is a 
smooth complete intersection
by a generalization of the theorem of Bertini
\cite[Theorem 1.10]{MR1201388}.

Since $X$ is a complete intersection,
the differential $ds$ of the section $s$ induces an isomorphism
\begin{align} \label{eq:normal}
  N_{X/F} \cong \cE|_X.
\end{align}
By taking the determinant of the exact sequence
\begin{align} \label{eq:tangent}
 0 \to T_X \to T_F|_X \to N_{X/F} \to 0,
\end{align}
one obtains an isomorphism
\begin{align}
 \det T_X \cong \det T_F|_X \otimes \det^{-1} \cE|_X.
\end{align}
Hence
$\det V \cong \det (\frakg/\frakp)$
is a sufficient condition for $\det T_X \cong \cO_X$,
which is necessary if the restriction map $\Pic F \to \Pic X$ is injective.

The exact sequence
\begin{align} \label{eq:omegaj}
 0 \to \Sym^j \cE^\dual|_X \to \cdots \to
 \Sym^{j-k} \cE^\dual \otimes \Omega^k_F |_X \to \cdots \to
 \Omega^j_F|_X \to \Omega^j_X \to 0
\end{align}
obtained as the $j$-th exterior power of the exact sequence
\begin{align} \label{eq:cotangent}
 0 \to \cE^\dual|_X \to \Omega^1_F|_X \to \Omega^1_X \to 0
\end{align}
dual to \pref{eq:tangent}
gives the spectral sequence
\begin{align} \label{eq:conormalss}
  E_1^{-q,p} = H^p \lb \Sym^q \cE^\vee \otimes \Omega_F^{j-q}|_X\rb 
    \Rightarrow H^{p-q} \lb \Omega_X^j \rb.
\end{align}

The Koszul resolution
\begin{align} \label{eq:Koszul}
 0
  \to \wedge^{\rank \cE} \cE^\dual
  \to \cdots
  \to \cE^\dual
  \to \cO_F
  \to \cO_X
  \to 0
\end{align}
gives the spectral sequence
\begin{align}
  \label{eq:koszulss}
  E_1^{-q,p} = H^p \lb \wedge^q \cE^\vee \otimes \cG \rb 
    \Rightarrow H^{p-q} \lb \cG|_X \rb
\end{align}
for any coherent sheaf $\cG$ on $F$.

Together with the Hodge symmetry
$h^{p,q}(X)=h^{q,p}(X)$
and
the obvious fact
that $h^{p,q}(X)=0$ unless $0\le p \le \dim X$,
the Hodge numbers are often determined only from dimensions of the cohomology groups 
on the $E_1$-page,
although there are cases where 
one should look at morphisms more carefully.

The topological Euler number $\chi(X)$
can be computed by
\begin{align} \label{eq:CGB}
\chi(X) = \int_X c(T_X) = \int_F \frac{c(T_F)}{c(\cE)}\ \ctop(\cE),
\end{align}
where $c(\cG)$ and $\ctop(\cG)$ denote the total and the top Chern classes.
The first equality in \pref{eq:CGB} is the Chern--Gauss--Bonnet theorem,
and the second equality comes from \pref{eq:normal} and \pref{eq:tangent}.
From the splitting principle, it follows
\begin{align}
  \label{eq:euler}
\chi(X) =\int_F \frac{\prod_{\mu\in\Delta\lb\frakg/\frakp^\dual\rb} \lb 1+c_1(\cL_\mu) \rb}
{\prod_{\nu\in\Delta(V^\dual)} \lb 1+ c_1(\cL_\nu) \rb} 
\prod_{\lambda\in\Delta(V^\dual)} c_1(\cL_\lambda),
\end{align}
where 
$\cL_\lambda = \cE_{V^B_\lambda}^\dual$ is a line bundle on $G/B$ for any 
weight $\lambda$,
and the integrand is an element of
$H^*(F,\bZ)$
considered as a subgroup of
$H^*(G/B,\bZ)$
by the pull-back
along the natural projection
$
G/B \to F;
$
an element in $H^*(F,\bZ)$ is described
as a polynomial in
$
x_i
\coloneqq c_1(\cL_{\omega_i})
\in H^2(G/B,\bZ)
$
for $i=1,\dots,r$.
Note that
$
c_1(\cL_\lambda) = \sum_{i=1}^r \lambda_i x_i.
$
One can perform the integral in terms of representation theory by using the following \pref{lm:monomial}.

\begin{lemma}[{cf.~e.g.~\cite[Lemma 6.3.2]{MR1038279}}]
  \label{lm:monomial}
For a monomial $x_i x_j \dots x_k\in H^{2l}(G/B,\bZ)$ of degree $l$ and 
a Schubert cycle $[X_w]\coloneqq [\overline{BwB/B}] \in H_{2l}(G/B,\bZ)$ 
associated with $w \in W$ of length $l$, 
one has
\begin{align}
  x_i x_j \dots x_k [X_w] = \sum
  \killing{\lambda_{i}}{\beta_{(1}^\dual}
  \killing{\lambda_{j}}{\beta_2^\dual} \cdots 
  \killing{\lambda_{k}}{\beta_{l)}^\dual},
\end{align}
where the sum runs over all collections $\beta_1, \dots, \beta_l \in \Delta^+$
such that
\begin{align}
  w = \sigma_{\beta_1}\sigma_{\beta_2}\dots \sigma_{\beta_l}
\end{align}
is a reduced expression,
and the parentheses on the right-hand side denote the symmetrized product.
\end{lemma}


\section{Rational homogeneous spaces of exceptional types}
 \label{sc:G2}
Let $G$ be the complex simple Lie group of exceptional type.
We call a rational homogeneous space $G/P$ an \emph{exceptional flag variety},
or an \emph{exceptional Grassmannian} if $P$ is maximal.

For type $G_2$,
there are three homogeneous spaces
$G/P_1$, $G/P_2$ and $G/B$
associated with the crossed Dynkin diagrams
\dynkin[parabolic=1]G2, \dynkin[parabolic=2]G2 and \dynkin[mark=x]G2 respectively.
The sets of weights of the irreducible representations
$V^{P}_{(a,b)}$
with the highest weight
$
 (a,b) \coloneqq a \omega_1 + b \omega_2
$
are given by
\begin{align*}
 \Delta \lb V^{P_1}_{(a,b)} \rb
  &= \lc (a+j, b-2j) \relmid j=0, 1, \dots , b \rc, \\
 \Delta \lb V^{P_2}_{(a,b)} \rb
  &= \lc (a-2j, b+3j) \relmid j=0, 1, \dots , a \rc, \\
 \Delta \lb V^{B}_{(a,b)} \rb
  &= \lc (a, b) \rc.
\end{align*}
In particular,
the dimensions and the determinants of these representations
are given as follows:
\begin{align}
\begin{array}{ccc}
\toprule
 \text{representation} & \text{dimension} & \text{determinant} \\
\midrule
 V^{P_1}_{(a,b)} & b+1 & (a(b+1) + b(b+1)/2, 0) \\
 V^{P_2}_{(a,b)} & a+1 & (0, (a+1)b+3a(a+1)/2) \\
 V^{B}_{(a,b)} & 1 & (a, b) \\
\bottomrule
\end{array}
\end{align}
By considering the action of the nilpotent parts on the roots, 
one can directly observe that
the representations $\lb \frakg/\frakp \rb^\vee$ are given by
\begin{align}
 \lb \frakg/\frakp_1 \rb^\vee &\cong V^{P_1}_{(-1,3)} + V^{P_1}_{(1,0)}, \\
 \lb \frakg/\frakp_2 \rb^\vee &\cong V^{P_2}_{(1,-1)} 
  + V^{P_2}_{(0,1)} + V^{P_2}_{(1,0)}, \\
 \lb \frakg/\frakb \rb^\vee &\cong V^{B}_{(2,-3)}\oplus V^{B}_{(-1,2)} 
  + V^{B}_{(1,-1)} + V^{B}_{(0,1)}
  + V^{B}_{(-1,3)} + V^{B}_{(1,0)},
\end{align}
which agree with the formula for the cotangent bundles in \pref{lem:cotangent}.
The determinants are given by
\begin{align}
 \det \lb \frakg/\frakp_1 \rb^\vee &\cong V^{P_1}_{(3,0)}, \\
 \det \lb \frakg/\frakp_2 \rb^\vee &\cong V^{P_2}_{(0,5)}, \\
 \det \lb \frakg/\frakb \rb^\vee  &\cong V^{B}_{(2,2)}.
\end{align}

\pref{th:classification} for $G_2$-Grassmannians is
an easy consequence of these facts.
For example,
to obtain a Calabi--Yau 3-fold
from $(G/P_1, \cE_V^\dual)$,
the completely reducible representation $V$ must satisfy
$\dim V = 2$ and $\det V = V^{P_1}_{(3,0)}$
since $\Pic G/P_1 \cong \bZ$ must inject to $\Pic X$.
If $V$ is decomposable, then
$V \cong V^{P_1}_{(1,0)} \oplus V^{P_1}_{(2,0)}$
is the only choice,
and if $V$ is indecomposable, then
$V \cong V^{P_1}_{(1,1)}$ is the only choice.

For the Calabi--Yau 3-fold $X$
associated with
$
\lb G/P_1,\cE_{(1,1)} \rb,
$
one can use the spectral sequences 
\pref{eq:extension},
\pref{eq:conormalss}, and
\pref{eq:koszulss}
to prove
$h^{0,1} = h^{0,2} = 0$,
$h^{1,1} = 1$, and
$h^{1,2} = 50$
by hand.
The Koszul resolution \pref{eq:Koszul} allows us to compute
the cohomology of $\cO_X(i)$,
which together with the Hirzebruch--Riemann--Roch theorem
\begin{align}
 \label{eq:HRR}
 \chi(\cO_X(i)) = \frac{1}{6} \deg X \cdot i^3 + \frac{1}{12} c_2(X) \cdot i
\end{align}
implies
$
\deg X = 42$
and
$
c_2(X) = 84.
$

Similar calculations
aided by a Mathematica package \cite{1912.10969}
give Theorems \ref{th:classification},
\ref{th:Fano}, and
\ref{th:classification2}.
The conditions
$\rank \cE = \dim G/P -3$ and $c_1(\cE) = c_1(G/P)$
are strong,
and
many of 439 exceptional flag varieties
are eliminated quickly.
The topological invariants are calculated
by using the spectral sequences
\pref{eq:extension},
\pref{eq:conormalss},
\pref{eq:koszulss},
and the Chern--Gauss--Bonnet theorem 
\pref{eq:euler} in one case.
All calculations are recorded
in the ancillary files to this paper.

\bibliographystyle{amsalpha}
\bibliography{bibs}
\end{document}